\documentclass[12pt,letterpaper]{article}
\usepackage{geometry}                % See geometry.pdf to learn the layout options. There are lots.
\usepackage{graphicx}
\usepackage{amssymb}
\usepackage{epstopdf}

\usepackage{amsmath}
\usepackage[colorlinks=true,linkcolor=red,citecolor=blue]{hyperref}
 \usepackage{amsmath}
 \usepackage[latin1]{inputenc}
 \usepackage[T1]{fontenc}
\usepackage{times}
\usepackage{epsfig}
 \usepackage{dsfont}
 \usepackage{natbib}
 \usepackage{mathrsfs}

\newtheorem{theorem}{Theorem}[section]

\newtheorem{corollary}[theorem]{Corollary}

\newtheorem{remark}[theorem]{Remark}

\numberwithin{equation}{section}

\title{A Strong Invariance Theorem of the Tail Empirical Copula Processes}
\author{Salim BOUZEBDA$^{1}$\footnote{e-mail:
salim.bouzebda@upmc.fr} \hbox{  } and \hbox{} Tarek ZARI$^{2}$\footnote{e-mail: zaritarek@gmail.com}\\
 $^{1}$L.S.T.A.,  Universit\'e Pierre et Marie Curie-Paris VI\\
    4 place Jussieu
    75252 Paris Cedex 05, France\\
$^{2}$ Universit\'e Hassan II - FSJES A\^{\i}n Sba\^{a}\\
 Beausite, B.P. 2634 A\^{\i}n Sba\^{a} - Casablanca\\
                   20000 Casablanca,
                   Maroc.}

\begin{document}
\maketitle
\begin{abstract}
\noindent We study the behavior of bivariate empirical copula process $\mathbb{G}_n(\cdot,\cdot)$ on pavements $[0,k_n/n]^2$ of $[0,1]^2,$ where $k_n$ is a sequence of positive constants fulfilling some conditions. We provide a upper bound for the strong approximation of $\mathbb{G}_n(\cdot,\cdot)$  by a Gaussian process when $k_n/n \searrow \gamma$ as $n\rightarrow \infty,$ where $0 \leq \gamma \leq 1.$\\
\noindent{\small {\bf AMS Subject Classifications}: Primary 60F17 ; secondary 62G20 ; 62H10 ; 	60F15.} \\
 \noindent{\small {\bf Keywords}: Empirical Copula processes ; Strong invariance principles ; Rates of convergences ; Gaussian processes. }
 \end{abstract}
 
\section{Introduction} Let $\{(X_n,Y_n):n\geq 1\}$ be independent replic{\ae} of a random vector $(X,Y)$ with distribution function [df] $\mathbb{F}(x,y)=\mathbb{P}(X\leq x, Y\leq y)$. We assume that the corresponding
marginal df's $G(x)=\mathbb{P}(X\leq x)$ and $H(y)=\mathbb{P}(Y\leq y)$ are continuous. It is well know that there exists a distribution function $\mathbb{C}(\cdot,\cdot)$ with uniform marginals on $[0,1]^2$ such that
\begin{equation*}
\mathbb{F}(x,y)= \mathbb{C}(G(x),H(y))~~~\hbox{for}~~ (x,y)\in
\mathbb{R}^2.
\end{equation*}
See \cite{Sklar1959,Sklar1973}, \cite{Moore1975}, \cite{deheuvels1979a}. The function $\mathbb{C}(\cdot,\cdot)$ is called the copula associated with $\mathbb{F}(\cdot,\cdot)$ (some authors called it the dependence function). This function fulfills the following identity
\begin{equation}\label{def copule}
\mathbb{C}(u,v)= \mathbb{F}(G^{-1}(u),H^{-1}(v))~~~\hbox{for}~~(u,v)\in [0,1]^2,
\end{equation}
where $G^{-1}(u)=\inf\{x: G(x)\geq u\}$ and $H^{-1}(v)=\inf\{y: H(y)\geq v\}$ are the quantile functions pertaining, respectively, to $G(\cdot)$ and $H(\cdot)$. In the monographs by
\cite{nelsen2006} and \cite{Joe1997} the reader may find detailed
ingredients of the modelling theory as well as surveys of the commonly
used copulas.\\
The empirical counterparts of $\mathbb{F}(\cdot,\cdot),G(\cdot)$ and $H(\cdot),$ based upon $(X_1,Y_1),\ldots,(X_n,Y_n),$ are given, respectively, for each $n\geq 1$ and $x,y \in \mathbb{R},$ by
%\begin{equation*}
$$\mathbb{F}_n(x,y):=\frac{1}{n}\sum_{i=1}^{n}{\rm 1\!I}\{X_{i}\leq
x, Y_{i}\leq y\}, $$
%\end{equation*}
\begin{equation*}
G_{n}(x):=\mathbb{F}_n(x,\infty)=\frac{1}{n}\sum_{i=1}^{n}{\rm 1\!I}\{X_{i}\leq
x\}~~\hbox{\rm{and}}~~ H_{n}(y):=
\mathbb{F}_n(\infty,y)=\frac{1}{n}\sum_{i=1}^{n}{\rm 1\!I}\{Y_{i}\leq
y\},
\end{equation*}
where ${\rm 1\!I}\{\cdots\}$ denotes the indicator function of the set $\{\cdots\}$.
In view of the characterization (\ref{def copule}), we define an empirical copula function of $\mathbb{F}_n(\cdot,\cdot)$ by
\begin{equation}\label{copule empirique}
 \mathbb{C}_{n}(u,v) := \mathbb{F}_{n}(G_{n}^{-1}(u),H_{n}^{-1}(v)),
\end{equation}
where $G_n^{-1}(u):=\inf\{x: G_n(x)\geq u\}$ and $H_n^{-1}(v):=\inf\{y: H_n(y)\geq v\},$ for $0\leq u,v\leq 1$ and $n\geq 1,$ denote the empirical quantile functions of  $G_n(\cdot)$ and $H_n(\cdot)$, respectively. \\
We may now define the empirical copula process $\mathbb{G}_n(\cdot,\cdot)$ by setting
\begin{equation}\label{processus de copule}
\mathbb{G}_n(u,v):= n^{1/2}(\mathbb{C}_n(u,v)-
\mathbb{C}(u,v)),~~~~\hbox{for}~~(u,v)\in [0,1]^2.
\end{equation}
The asymptotic behavior of the copula process $\mathbb{G}_n(\cdot,\cdot)$ has been investigated in various setting by many authors in the context of process convergence. \cite{deheuvels1979a} investigated the
consistency of $\mathbb{C}_{n}(\cdot,\cdot)$ and \cite{Deheuvels1981b,deheuvels1981} proved a functional central limit theorem for $\mathbb{G}_n(\cdot,\cdot)$ in the particular case of independent  margins. \cite{Ruch1,Rush2}, \cite{stute1987}
proved weak convergence of the empirical copula process $\mathbb{G}_n(\cdot,\cdot)$ in the Skorokhod space $D([0,1]^2).$ \cite{Wellner1996} established weak convergence in the space $\ell^{\infty}([a,b]^2),$ when $0< a< b<1.$ \cite{fermanianradulovicdragan2004} showed  that the weak convergence of $\mathbb{G}_n(\cdot,\cdot)$ to a centered Gaussian process $\mathbb{G}(\cdot,\cdot)$ holds on $\ell^{\infty}([0,1]^2),$ when $\mathbb{C}(\cdot,\cdot)$ has continuous partial derivatives on $[0,1]^2.$ In particular, these conditions are satisfied under the independence assumption of margins, i.e., for $(u,v)\in [0,1]^{2}$
$$
\mathbb{C}(u,v)=uv.
$$
A natural question arises then: which is the rate for the strong approximation, $V_n$ say, of empirical copula processes $\mathbb{G}_n(\cdot,\cdot)$ by sequences of Gaussian processes $\mathbb{B}_n(\cdot,\cdot)$ such that
\begin{equation*}
 \sup_{(u,v)\in [0,1]^2} \mid \mathbb{G}_n(u,v)-\mathbb{B}_n(u,v) \mid = \mathcal{O}(V_n),~~~~\mbox{a.~s.}?
\end{equation*}
This is known in the scientific literature under the name of \emph{invariance principle}.\\ Under the independence assumption of margins, \cite{deheuvels2006} proved that the strong invariance principle holds with $V_n= n^{-1/4} (\log n)^{1/2} (\log \log n)^{1/4},$ where $\mathbb{B}_n(\cdot,\cdot)$ is a sequence of bivariate \emph{tied-down Brownian bridges }(see (\ref{tied down}) below). In the multivariate case, i.e., the dimension $d\geq 3$, \cite{Deheuvelsref2009} showed  that the rate for the strong approximations of empirical copula process is  $V_n=n^{-1/(2d)}(\log n)^{2/d}$ where $\mathbb{B}_n(\cdot)$ is  a sequence of $d$-variate \emph{ tied-down Brownian bridges} (see, Theorem 2.1, p.140 of \cite{Deheuvelsref2009}).\\
\noindent Set $\log_1 u=\log_{1}u=\log(u\vee e)$, and $\log_p u= \log_{1}(\log_{p-1}u)$ for $p \geq 2$.
Throughout this paper, $\{k_n\}_{n=1}^{\infty}$ will denote a sequence of positive constants such that for all integers $n>1,$
\begin{description}
\item(H.1)  $0<k_n \leq n$;
\item(H.2)  $k_n \uparrow$~~as~~$n\rightarrow\infty$;
\item(H.3)  $k_n/n \searrow \gamma $~~with~~$0 \leq \gamma \leq 1$~~~as~~$n\rightarrow\infty$;
\item(H.4) $ k_n/\log_2 n \rightarrow\infty $,~~as~~$n\rightarrow\infty$.
\end{description}
In the present paper, we  are mainly concerned with the empirical copula process defined, in terms of a sequence $\{k_n\}_{n=1}^{\infty}$, for each $n\geq 1$, by
\begin{eqnarray} \label{processus G_n sur le pave}
\mathbb{G}_n^{*}(u,v):=
\mathbb{G}_n\left(u\frac{k_n}{n},v\frac{k_n}{n}\right)~~~~\hbox{for}~~(u,v)\in[0,1]^{2}.
\end{eqnarray}
The remainder of the paper is organized as follows. In \S \ref{section 2}, we state the main result
concerning the rate of uniform almost sure convergence of the process $\{\mathbb{G}_n^{*}(u,v): (u,v)\in [0,1]^{2}\}$, $n\geq 1$, defined in (\ref{processus G_n sur le pave}), in the case of independent and general marginals. These results will be used, in \S\ref{testing}, to derive some asymptotic properties of the statistic considered in this paper for testing tail independence. We also indicate the
basic technical tools needed for establishing this result. In \S\ref{section3}, we apply our result to the smoothed local empirical copula process, where we limit ourselves to the case of independent marginals.
\section{Behavior of empirical copula process on pavements}\label{section 2}
\subsection{Gaussian Process}

Our main aim is to provide  a strong approximation of the process $\mathbb{G}_n^{*}(\cdot,\cdot)$ based upon $$(X_{1},Y_{1}),\ldots,(X_{n},Y_{n})$$ by a sequence of Gaussian processes. To award this goal, we first need to introduce several approximating Gaussian processes, we refer for example to  \cite{Deheuvels2007}. \\
Let $\{\mathbf{W}(s,t): s\geq 0, ~t\geq 0 \}$ be a two-time parameter Wiener process (or \emph{Brownian sheet}), namely, a centered Gaussian process, with continuous trajectories and covariance function
given by
$$ \mathbb{E}(\mathbf{W}(s_1,t_1)\mathbf{W}(s_2,t_2))= (s_1\wedge s_2) (t_1\wedge t_2),~~~~\hbox{for}~~s_1,s_2,t_1,t_2\geq 0,      $$
where $(s_1\wedge s_2)=\min(s_1,s_2).$\\
A bivariate \emph{ Brownian bridge} is defined, in terms of the two-time Wiener process $\mathbf{W}(\cdot,\cdot),$ via
\begin{equation}
\mathbf{B}(s,t)=\mathbf{W}(s,t)-st\mathbf{W}(1,1),~~~~\hbox{for}~~(s,t)\in [0,1]^2.
\end{equation}
This process has continuous sample paths and fulfills
$$ \mathbb{E}(\mathbf{B}(s,t))=0,~~ \mathbb{E}(\mathbf{B}(s_1,t_1)\mathbf{B}(s_2,t_2))= (s_1\wedge s_2) (t_1\wedge t_2)- \prod_{i=1}^2 s_it_i,~~~~\hbox{for}~~(s_1,s_2,t_1,t_2)\in [0,1]^4.      $$
A bivariate \emph{tied-down Brownian bridge} is defined, in terms of the bivariate Brownian bridge $\mathbf{B}(\cdot,\cdot),$ via
\begin{equation}\label{tied down}
\mathbb{B}(s,t)=\mathbf{B}(s,t)-s\mathbf{B}(1,t)-t\mathbf{B}(s,1),~~~~\hbox{for}~~(s,t)\in [0,1]^2.
\end{equation}
It is a centered process with continuous sample paths and  covariance function given by
$$  \mathbb{E}(\mathbb{B}(s_1,t_1)\mathbb{B}(s_2,t_2))= \{s_1\wedge s_2-s_1s_2\} \{t_1\wedge t_2- t_1t_2 \},~~~~\hbox{for}~~(s_1,s_2,t_1,t_2)\in [0,1]^4.      $$
\subsection{Asymptotic Theory}
Throughout the sequel, we assume that $X$ and $Y$ are mutually independent with continuous distribution functions $G(\cdot)$ and $H(\cdot).$ Thus,  $\{U_i= G(X_i)\}_{1\leq i\leq n}$ and $\{V_i= H(Y_i)\}_{1\leq i\leq n}$ are two independent sequences of  independent and identically distributed uniform $(0,1)$ random variables. For all $(u,v)\in
[0,1]^2$, the distribution function $\mathbb{T}(\cdot,\cdot)$ associated with $(U,V)$ fulfills the following identity
$$ \mathbb{T}(u,v) := \mathbb{P}(U \leq u, V\leq v)= \mathbb{F}(G^{-1}(u),H^{-1}(v)) = \mathbb{C}(u,v)=uv. $$
We define, for each $n\geq1$ and $0\leq u,v \leq 1,$ the
empirical counterparts of $\mathbb{T}(\cdot,\cdot)$ and the empirical marginals based on $\{(U_i,V_i)\}_{1 \leq i \leq n},$ respectively, by setting
\begin{eqnarray}
 \mathbb{T}_n(u,v)&:=& \frac{1}{n} \sum_{i=1}^{n} {\rm 1\!I}\{U_i
\leq u, V_i \leq v\}:= \mathbb{F}_n(G^{-1}(u),H^{-1}(v)),\\
\mathbb{U}_n(u)&=& \mathbb{T}_n(u,1)=G_n(G^{-1}(u)),\\
\mathbb{V}_n(v)&=& \mathbb{T}_n(1,v)=H_n(H^{-1}(v)).
\end{eqnarray}
The empirical quantile functions of $\mathbb{U}_n(\cdot)$ and $\mathbb{V}_n(\cdot)$ are given, for $n\geq 1$ and $0\leq u,v \leq 1,$ by
\begin{eqnarray} \label{quantile U}
\mathbb{U}_n^{-1}(u) &=& \inf \{s \geq 0: \mathbb{U}_n(s)\geq u \}
= G(G_n^{-1}(u)),\\
\label{quantile V} \mathbb{V}_n^{-1}(v) &=& \inf \{t \geq 0:
\mathbb{V}_n(t)\geq v \} = H(H_n^{-1}(v)).
\end{eqnarray}
 Consider the empirical processes defined, respectively, for each $n\geq 1$ and $0\leq u,v \leq 1,$ by \begin{eqnarray}
 \alpha_n(u,v)&=& n^{1/2}(\mathbb{T}_n(u,v)-uv),\\
\label{processus U_n} \alpha_{n;\mathbb{U}}(u)&=&\alpha_n(u,1)=
n^{1/2}(\mathbb{U}_n(u)-u),\\
\alpha_{n;\mathbb{V}}(v)&=&\alpha_n(1,v)=
n^{1/2}(\mathbb{V}_n(v)-v),\\
\beta_{n;\mathbb{U}}(u)&=& n^{1/2}(\mathbb{U}_n^{-1}(u)-u
),\\
\label{processus beta_n de V} \beta_{n;\mathbb{V}}(v)&=&
n^{1/2}(\mathbb{V}_n^{-1}(v)-v ).
\end{eqnarray}
In view of the definition (\ref{copule empirique}) of $\mathbb{C}_n(\cdot,\cdot),$ the relation between $\mathbb{C}_n(\cdot,\cdot)$ and $\mathbb{T}_n(\cdot,\cdot)$ is given, for each $n\geq 1$ and   $0\leq
u,v\leq 1,$ by
\begin{eqnarray}
\mathbb{T}_n(\mathbb{U}_n^{-1}(u),\mathbb{V}_n^{-1}(v))=  \mathbb{F}_n(G_n^{-1}(u),H_n^{-1}(v))= \mathbb{C}_n(u,v).         \end{eqnarray}
In order to study the process $\mathbb{G}_n^{*}(\cdot,\cdot)$ defined in (\ref{processus G_n sur le pave}), and in view of (\ref{processus U_n})-(\ref{processus beta_n de V}), we define the tail empirical and quantile processes in terms of the sequence $\{k_n\}_{n=1}^{\infty},$ for each $n\geq 1$ and   $0\leq
u,v\leq 1,$ by
\begin{eqnarray}
 \alpha_n^{*}(u,v)&:=& \alpha_n\left(u\frac{k_n}{n},v\frac{k_n}{n}\right),\\
 \alpha_{n;\mathbb{U}}^{*}(u) &:=& \left(\frac{k_n}{n}\right)^{-1/2}
\alpha_{n;\mathbb{U}}\left(u\frac{k_n}{n}\right),\\
 \alpha_{n;\mathbb{V}}^{*}(v) &:=& \left(\frac{k_n}{n}\right)^{-1/2}
\alpha_{n;\mathbb{V}}\left(\frac{v k_n}{n}\right),\\
\label{beta*u} \beta_{n;\mathbb{U}}^{*}(u)&:=&
\left(\frac{k_n}{n}\right)^{-1/2}
\beta_{n;\mathbb{U}}\left(u\frac{k_n}{n}\right),\\
\label{beta*v} \beta_{n;\mathbb{V}}^{*}(v)&:=&
\left(\frac{k_n}{n}\right)^{-1/2} \beta_{n;\mathbb{V}}\left(\frac{v
k_n}{n}\right).
\end{eqnarray}
The tail of the empirical process of $\alpha_{n;\bullet}(\cdot)$ and the tail of the quantile  process $\beta_{n;\bullet}(\cdot)$ play  major and important role in statistics, for instance,  the nonparametric statistics.  This importance explains the huge variety of existing results in this fields, we may refer to \cite{sDeheuvels1997}, \cite{masoneinmahl1988, masonein1988} and the references therein. The strong approximations of  the processes  $\alpha_n^{*}(\cdot,\cdot),\alpha_{n;\bullet}^{*}(\cdot)$ and $\beta_{n;\bullet}^{*}(\cdot)$  are given in \cite{Mason1987}, \cite{Mason1988}, \cite{Csorgo1988}.\\
Keep in mind  the definition of $\mathbb{G}_{n}^{*}(\cdot,\cdot)$ in (\ref{processus G_n sur le pave}), for each $n\geq 1$ and $0\leq u,v \leq 1,$ we have
\begin{eqnarray}
\nonumber \mathbb{G}_{n}^{*}(u,v)&=& \mathbb{G}_{n}\left(\frac{u
k_n}{n},\frac{v k_n}{n}\right)\\
\nonumber &=&  n^{1/2}\left\{\mathbb{C}_n\left(\frac{u
k_n}{n},\frac{v k_n}{n}\right)-\frac{u k_n}{n}\frac{vk_n}{n} \right\}.
\end{eqnarray}
This process can be decomposed as follows
\begin{eqnarray}
\nonumber
\lefteqn{\mathbb{G}_{n}^{*}(u,v)}\\
\nonumber&=& n^{1/2}\left\{\mathbb{T}_n\left(\mathbb{U}_n^{-1}\left(\frac{u k_n}{n}\right),\mathbb{V}_n^{-1}\left(\frac{vk_n}{n}\right)\right)-\frac{u k_n}{n}\frac{vk_n}{n} \right\}\\
\nonumber &=& \alpha_n\left(\mathbb{U}_n^{-1}\left(\frac{u
k_n}{n}\right),\mathbb{V}_n^{-1}\left(\frac{vk_n}{n}\right)\right) + n^{1/2}\left[\mathbb{U}_n^{-1}\left(\frac{u k_n}{n}\right)\mathbb{V}_n^{-1}\left(\frac{vk_n}{n}\right)-\frac{u k_n}{n}\frac{vk_n}{n}\right]\\
\nonumber &=& \alpha_n\left(\mathbb{U}_n^{-1}\left(\frac{u
k_n}{n}\right),\mathbb{V}_n^{-1}\left(\frac{vk_n}{n}\right)\right)+u
\frac{k_n}{n} \beta_{n;\mathbb{V}}\left(v\frac{k_n}{n}\right)+v\frac{k_n}{n}\beta_{n;\mathbb{U}}\left(u\frac{k_n}{n}\right)\\
&&+n^{-1/2}\beta_{n;\mathbb{U}}\left(u\frac{k_n}{n}\right)\beta_{n;\mathbb{V}}\left(v\frac{k_n}{n}\right).
\end{eqnarray}
Our approximations will be based on the following fact, which combines results of \cite{masoneinmahl1988} and \cite{masonein1988}. For convenience, we will denote \emph{sup-norm} of a bounded function $f(\cdot),$ defined on $I=[0,1]$ or $I=[0,1]^2,$ by $\| f \|= \sup_{x\in I}| f(x)|.$ The next fact, due to \cite{masoneinmahl1988} provides a law of the iterated logarithm of the local quantile process.\\
{\bf Fact 1.}
Under (H.1)-(H.4), we have, with probability 1,
\begin{equation}
\begin{array}{lcl}
 \limsup_{n\rightarrow\infty}  (\log_2 n)^{-1/2} \|\beta_{n;\bullet}^*\| &=&
2^{1/2} (1-\gamma)^{1/2}~~~\mbox{if}~~0 \leq \gamma \leq 1/2,  \\
&&\\
\label{L.I.L tail quantile1}
&=& 2^{-1/2} \gamma^{-1/2}~~~~~~~\mbox{if}~~1/2 < \gamma \leq
1 .
\end{array}
 \end{equation}
\cite{masonein1988} established Bahadur-Kiefer type representation, which relates
to the sum of the local uniform empirical process and the local quantile process. Let
\begin{eqnarray} \label{expression de R_i}
 {\bf{R}}_{n;\bullet}(k_{n})&:=& \sup_{s \in
[0,\frac{k_{n}}{n}]} |\alpha_{n;\bullet}(s)+\beta_{n;\bullet}(s)|,\\
 \label{expression de r_i}
 r_{n}&:=& n^{-1/2} k_{n}^{1/4} ( \log_2 n)^{1/4}
\left( \log_{1} (k_{n}) + 2 \log_2 n\right)^{1/2}.
\end{eqnarray}
In the sequel, we need the following fact due to \cite{masonein1988}.\\
{\bf Fact 2.}
\label{kiefer local}
Let $\{k_n\}_{n=1}^{\infty}$ a sequence of positive constants which satisfy the assumptions
(H.1)-(H.4).
\begin{description}
    \item(i) When $\gamma = 0$,  we have, with probability 1,
\begin{gather} \label{gamma 0}
 \limsup_{n\rightarrow \infty} r_{n}^{-1} {\bf{R}}_{n;\bullet}(k_n)\leq
2^{1/4} .
\end{gather}
 In addition, when $\log_{1}(k_n) / \log_2 n \rightarrow
\infty~~\mbox{as}~~n\rightarrow\infty$, we obtain an equality on (\ref{gamma
0}).
    \item(ii) When $0< \gamma \leq 1$, we have, with probability 1,
\begin{eqnarray*}\label{kiefer local2}
\limsup_{n\rightarrow \infty}r_{n}^{-1} {\bf{R}}_{n;\bullet}(k_n) &=&
2^{1/4} (1-\gamma)^{1/4},~~~~0< \gamma \leq 1/2, \\
&=& 2^{-1/4} \gamma^{-1/4},~~~~\hspace{.8cm}1/2< \gamma \leq 1.
\end{eqnarray*}
\end{description}
In view of  (\ref{quantile U}), (\ref{quantile V}), (\ref{beta*u}) and
 (\ref{beta*v}), we have
\begin{eqnarray}
\label{tail quantile process1} \mathbb{U}_n^{-1}(u k_n/n)&:=& u
k_n/n + n^{-1} k_n^{1/2}\beta_{n;\mathbb{U}}^{*}(u),\\ \label{tail quantile
process2} \mathbb{V}_n^{-1}(v k_n/n)&:=& v k_n/n +
n^{-1}k_n^{1/2}\beta_{n;\mathbb{V}}^{*}(v).
\end{eqnarray}
Consider now the modulus of continuity $w_{n}(\cdot)$ of $\alpha_{n}(\cdot,\cdot)$ defined by
\begin{equation}\label{definition du module d oscillation}
 w_{n}(h_{n}) := \sup_{L \in \mathcal{R} : \mid L \mid \leq h_{n}} \mid
\alpha_{n}(L)\mid ~~~\hbox{for}~~h_{n} \in (0,1), \end{equation}
where
\begin{eqnarray*}
  \mathcal{R}&:=&
\left\{[{\bf s},{\bf t}]= [s_1,t_1]\times[s_2,t_2]: 0
\leq s_i \leq t_i \leq 1~~\mbox{for}~~ i=1,2\right\},\\
|L|&=& |{\bf t}-{\bf s}|= \prod_{i=1}^2 |t_i-s_i|
\end{eqnarray*}
 and $h_{n}$ denotes a sequence of positive constants fulfilling the conditions of the following fact  due to \cite{sEinmahl1987}.\\
{\bf Fact 3.}
Let $\{h_n\}_{n=1}^{\infty}$ be a sequence of positive numbers on
$(0,1)$ with $h_n \downarrow 0$ as  $n\rightarrow
\infty$, such that
$$ i)  nh_n\uparrow \infty,~~~~ii) nh_n/ \log_{1}n \rightarrow \infty
,~~~~ iii) \log_{1}(1/h_n)/ \log_{2}n   \rightarrow \infty. $$
Then, with probability 1,
\begin{equation}  \label{oscillation einmahl}
\lim_{n\rightarrow
\infty} (2 h_n \log_{1}(1/ h_n))^{-1/2}w_{n}(h_n) = 1.
\end{equation}
 The proof of our result relies on the following oscillation inequality for bivariate empirical process,
which is mentioned in \cite{deheuvels2006}.\\
{\bf Fact 4.} \label{oscillation deheuvels}
For $0\leq u_1,v_1,u_2,v_2 \leq 1$, we have
\begin{eqnarray} \label{inegalite deheuvels}
\mid \alpha_n(u_1,v_1) - \alpha_n(u_2,v_2) \mid \leq 3 \times
w_{n}(|u_1 -u_2|\vee |v_1 -v_2|) .
 \end{eqnarray}
 For  each $n\geq 1$ and
$0\leq u,v \leq 1,$ set
\begin{eqnarray}
\nonumber
\alpha_{n;0}^{*}(u,v)&:=& \alpha_{n;0}\left(\frac{uk_n}{n},\frac{vk_n}{n} \right)\\
&:=& \alpha_n\left(\frac{uk_n}{n},\frac{vk_n}{n}\right)-
u\frac{k_n}{n}
\alpha_{n,\mathbb{V}}\left(v\frac{k_n}{n}\right)-v\frac{k_n}{n}\alpha_{n,\mathbb{U}}\left(u\frac{k_n}{n}\right).
\end{eqnarray}
We are now in position to study the behavior of the empirical copula process on  $[0,\frac{k_n}{n}]\times
[0,\frac{k_n}{n}].$ Our main result is summarized in the following theorem.

\begin{theorem}\label{theoreme pave}
Let $\left\{k_n\right\}_{n\geq1}$ be a sequence of positive numbers fulfilling the assumptions
{\rm(H.1)-(H.4)}. We have almost surely,
\begin{description}
    \item(i) when $ 0 < \gamma \leq 1/2$,
\begin{eqnarray} \label{resultat(i)}
\nonumber&& \limsup_{n\rightarrow\infty} n^{1/2} k_n^{-1/4}(\log_2 n)^{-1/4}
(\log_{1} n)^{-1/2} \|\mathbb{G}_{n}^{*}-\alpha_{n;0}^{*} \| \\&& \qquad \qquad\leq
[3\times 2^{-1/4}+ \gamma 2^{5/4}](1-\gamma)^{1/4},
\end{eqnarray}
   \item(ii) when $ 1/2 < \gamma \leq 1$,
  \begin{eqnarray}\label{resultat(ii)}
\nonumber &&\limsup_{n\rightarrow\infty} n^{1/2} k_n^{-1/4}(\log_2 n)^{-1/4}
(\log_{1} n)^{-1/2} \|\mathbb{G}_{n}^{*}-\alpha_{n;0}^{*} \| \\&&\qquad \qquad\leq
[3\times 2^{-3/4}+ \gamma 2^{3/4}] \gamma^{-1/4}.
\end{eqnarray}
\end{description}
\end{theorem}
\subsection*{Proof.} For each $n\geq 1$ and $0 \leq u,v \leq 1$,
we have
\begin{eqnarray}
 \mathbb{G}_{n}^{*}(u,v)-\alpha_{n;0}^{*}(u,v)
\nonumber & =&\left[\alpha_n\left(\mathbb{U}_n^{-1}\left(\frac{u
k_n}{n}\right),\mathbb{V}_n^{-1}\left(\frac{vk_n}{n}\right)\right)-\alpha_n \left(u\frac{k_n}{n},v\frac{k_n}{n} \right) \right] \\
\nonumber &  &              +
v\frac{k_n}{n}\left[\beta_{n;\mathbb{U}}\left(u\frac{k_n}{n}\right)+\alpha_{n;\mathbb{U}}\left(u\frac{k_n}{n}\right)\right]\\&&\nonumber
                 +u \frac{k_n}{n} \left[\beta_{n;\mathbb{V}}\left(v\frac{k_n}{n}\right)
                 +\alpha_{n;\mathbb{V}}\left(v\frac{k_n}{n}\right)\right] \\
\nonumber && + n^{-1/2}\beta_{n;\mathbb{U}}\left(u\frac{k_n}{n}\right)\beta_{n;\mathbb{V}}\left(v\frac{k_n}{n}\right)\\
\nonumber &
   = &R_{n;0}(u,v) + v\frac{k_n}{n} R_{n;\mathbb{U}}(u)
    + u\frac{k_n}{n} R_{n;\mathbb{V}}(v)+ R_{n}(u,v).
\end{eqnarray}
We study each quantity separately.  For the particular choice of
 $u_1=\mathbb{U}_n^{-1}\left(\frac{u k_n}{n}\right),$~$u_2=
u\frac{k_n}{n},~v_1= \mathbb{V}_n^{-1}\left(\frac{vk_n}{n}\right)$ and
$v_2=v\frac{k_n}{n}$ in (\ref{inegalite deheuvels}),
combined with (\ref{tail quantile process1}) and (\ref{tail quantile
process2}), we get that, with probability $1$, for  $n$ sufficiently large,
\begin{equation*}
\|R_{n;0}\| \leq 3 \times w_{n}(|u_1 -u_2|\vee |v_1 -v_2|).
   \end{equation*}
Observe that
$$|u_1 -u_2|=
n^{-1} k_n^{1/2} \|\beta_{n;\mathbb{U}}^{*}\| ~~~\hbox{et}~~~ |v_1 -v_2|=
n^{-1} k_n^{1/2} \|\beta_{n;\mathbb{V}}^{*}\| .
$$
First, we consider the case when  $ 0 \leq \gamma \leq 1/2.$ From the Fact 1, we infer that, almost surely, for $n$ sufficiently large,
$$\|\beta_{n;\mathbb{U}}^{*}\|=\|\beta_{n;\mathbb{V}}^{*}\|= \mathcal{O}\left( (\log_2
n)^{1/2}\right),$$ hence,  $$ |u_1 -u_2|\vee |v_1 -v_2|=
\mathcal{O}\left( n^{-1} k_n^{1/2}(\log_2 n)^{1/2} \right). $$
 Fix any $\epsilon > 0,$ and set $$h_n= (1+\epsilon) n^{-1} [2
(1-\gamma)]^{1/2} k_n^{1/2} (\log_2 n)^{1/2}.$$ By
 (\ref{oscillation einmahl}), we have almost surely,
\begin{eqnarray*}
&&\limsup_{n\rightarrow \infty} n^{1/2} k_n^{-1/4}(\log_2 n)^{-1/4} (\log_{1}
n)^{-1/2} w_{n}(|u_1 -u_2|\vee |v_1 -v_2|) \\&&\qquad\qquad\qquad= 2^{-1/4}
(1-\gamma)^{1/4} \sqrt{1+\epsilon}.
\end{eqnarray*}
Consequently, we have almost surely,
\begin{eqnarray}\label{R0}
\limsup_{n\rightarrow \infty} n^{1/2} k_n^{-1/4}(\log_2 n)^{-1/4}
(\log_{1} n)^{-1/2} \|R_{n;0}\| \leq 3 \times
 \left[\frac{1-\gamma}{2}\right]^{1/4}. \end{eqnarray}
Thought the sequel, set
$$ {\bf{V}}_n :=  n^{1/2} k_n^{-1/4}(\log_2 n)^{-1/4} (\log_{1}
n)^{-1/2}. $$
To study the second term $R_{n;\mathbb{U}}(\cdot),$ let us recall from
(\ref{expression de R_i}), (\ref{expression de r_i}) the definitions of $R_{n;.}(\cdot)$ and $r_{n}.$ By Fact 2, we have, almost surely,
\begin{eqnarray*}
\lefteqn{\limsup_{n\rightarrow \infty}\left\{ {\bf{V}}_n \|R_{n;\mathbb{U}}\|\right\}}\\ & \leq& \limsup_{n\rightarrow \infty}
\frac{k_n}{n}
 \times  \left[
\frac{\log_{1}( k_n)+ 2 \log_2 n}{\log_{1} n}\right]^{1/2} r_{n}^{-1}
{\bf{R}}_{n;\mathbb{U}}(k_n).
\end{eqnarray*}
Under (H.1)-(H.4), and keeping in mind  that the condition (H.3) implies that
$\left[ \frac{\log_{1} ( k_n)+ 2 \log_2 n}{\log_{1} n}\right]$ converges to
$1$ as $n\rightarrow\infty$, one can see  that, almost surely,
\begin{eqnarray}\label{R(n,U)}
\limsup_{n\rightarrow \infty} \left\{{\bf{V}}_n \|R_{n;\mathbb{U}}\| \right\} \leq 2^{1/4} \gamma
(1-\gamma)^{1/4}.
\end{eqnarray}
Similarly, using the same preceding arguments, we obtain, almost surely,
\begin{eqnarray}\label{R(n,V)}
\limsup_{n\rightarrow \infty} \left\{{\bf{V}}_n \|R_{n;\mathbb{V}}\|\right\} \leq 2^{1/4} \gamma (1-\gamma)^{1/4}.
\end{eqnarray}
Since $$\|\beta_{n;\mathbb{U}}\|=\|\beta_{n;\mathbb{V}}\|= \mathcal{O}\left( (k_n
\log_2 n)^{1/2}\right),$$ we then conclude that, almost surely,
\begin{eqnarray}\label{R(U,V)}
 \limsup_{n\rightarrow \infty} \left\{{\bf{V}}_n \|R_{n}\|\right\} = 0 .
\end{eqnarray}
By combining (\ref{R0}), (\ref{R(n,U)}), (\ref{R(n,V)}) with (\ref{R(U,V)}) we obtain
(\ref{resultat(i)}).\\
In the case when $1/2 < \gamma \leq 1 $, by using the same arguments to proof (\ref{resultat(i)}) and choosing
$$ h_n= (1+\epsilon) n^{-1}2^{-1/2} \gamma
^{-1/2} k_n^{1/2}(\log_2 n)^{1/2} \log_{1} n,$$
 to obain (\ref{resultat(ii)}). The proof of the Theorem is now completed. $\hfill \Box$\\

\noindent Recall the following result due to \cite{castellebonvalot1998}.
\begin{theorem}\label{bestapproximation}
On a suitable probability space, one can construct a sequence $\{
\mathbf{B}_{n}(u,v): (u,v) \in [0,1]^2 \}$ of
copula Brownian Bridges such that
 we may define the bivariate empirical process $\{
\alpha_n(u,v): (u,v) \in [0,1]^2 \}$, such that, for all positive $x$ and all $a$, $b\in [0,1]$,
the following inequality holds
\begin{eqnarray*}
&&\mathbb{P} \left(\sup_{0\leq u\leq a,0\leq v\leq b}|\alpha_{n}(u,v)-\mathbf{B}_{n}(u,v)|\geq n^{-1/2}(x+\Lambda_{1}\log(nab))\log(nab)\right)\\&&\qquad\qquad\leq \Lambda_{2}\exp(-\Lambda_{3}x),
\end{eqnarray*}
where $\Lambda_{1}$, $\Lambda_{2}$ and $\Lambda_{3}$ are absolute constants.
\end{theorem}
We define a sequence of \emph{tied-down Brownian bridges}, with the same law of $\left\{\mathbb{B}(s,t): (s,t) \in [0,1]^{2}\right\}$, by
  setting, for $n=1,2,\ldots,$ and $(s,t)\in [0,1]^2 ,$
\begin{equation}\label{GaussIndependen}
\mathbb{B}_{n}(s,t)=\mathbf{B}_{n}(s,t)-s\mathbf{B}_{n}(1,t)-t\mathbf{B}_{n}(s,1).  
\end{equation}
By combining Theorems \ref{theoreme pave} and \ref{bestapproximation}, we obtain the following corollary.
\begin{corollary}\label{COOOOO} Assume that the conditions of Theorem \ref{theoreme pave} and Theorem \ref{bestapproximation} hold. We have, almost surely,
\begin{description}
    \item(i)  when $ 0 < \gamma \leq 1/2$,
\begin{eqnarray} \label{resultat(i)}
\nonumber&&\limsup_{n\rightarrow\infty} n^{1/2} k_n^{-1/4}(\log_2 n)^{-1/4}
(\log_{1} n)^{-1/2} \|\mathbb{G}_{n}^{*}-\mathbb{B}_{n} \| \\&&\qquad\qquad \leq
[3\times 2^{-1/4}+ \gamma 2^{5/4}](1-\gamma)^{1/4},
\end{eqnarray}
   \item(ii) when $ 1/2 < \gamma \leq 1$,
  \begin{eqnarray}\label{resultat(ii)}
\nonumber&&\limsup_{n\rightarrow\infty} n^{1/2} k_n^{-1/4}(\log_2 n)^{-1/4}
(\log_{1} n)^{-1/2} \|\mathbb{G}_{n}^{*}-\mathbb{B}_{n} \|\\&&\qquad \qquad \leq
[3\times 2^{-3/4}+ \gamma 2^{3/4}] \gamma^{-1/4}.
\end{eqnarray}
\end{description}
\end{corollary}
\subsection{General case}
In this subsection we discuss the case when $X$ and $Y$ are dependent, i.e., for $(u,v)\in [0,1]^{2}$,
$$
\mathbb{C}(u,v)\neq uv.
$$
Assume that $\mathbb{C}(\cdot,\cdot)$ is twice continuously differentiable on $(0, 1)^2$ and all the
partial derivatives of second order are continuous on $[0,1]^2$. Then, By applying a Taylor series expansion, we have
\begin{eqnarray*}
\nonumber
\lefteqn{\mathbb{G}_{n}^{**}(u,v)}\\&=& n^{1/2}\left\{\mathbb{T}_n\left(\mathbb{U}_n^{-1}\left(\frac{u k_n}{n}\right),\mathbb{V}_n^{-1}\left(\frac{vk_n}{n}\right)\right)-\mathbb{C}\left(\frac{u k_n}{n},\frac{vk_n}{n}\right) \right\}\\
\nonumber &=& \widetilde{\alpha}_n\left(\mathbb{U}_n^{-1}\left(\frac{u
k_n}{n}\right),\mathbb{V}_n^{-1}\left(\frac{vk_n}{n}\right)\right) + n^{1/2}\left[\mathbb{C}\left(\mathbb{U}_n^{-1}\left(\frac{u k_n}{n}\right),\mathbb{V}_n^{-1}\left(\frac{vk_n}{n}\right)\right)-\mathbb{C}\left(\frac{u k_n}{n},\frac{vk_n}{n}\right)\right]\\
\nonumber&=&\widetilde{\alpha}_n\left(\mathbb{U}_n^{-1}\left(\frac{u
k_n}{n}\right),\mathbb{V}_n^{-1}\left(\frac{vk_n}{n}\right)\right) + n^{1/2}\left\{\mathbb{U}_n^{-1}\left(\frac{u k_n}{n}\right)-\frac{u k_n}{n}\right\}\frac{\partial \mathbb{C}}{\partial u}\left(\frac{u k_n}{n},\frac{vk_n}{n}\right)\\
&&+ n^{1/2}\left\{\mathbb{V}_n^{-1}\left(\frac{v k_n}{n}\right)-\frac{v k_n}{n}\right\}\frac{\partial \mathbb{C}}{\partial v}\left(\frac{u k_n}{n},\frac{vk_n}{n}\right)\\
%\end{eqnarray*}
%\begin{eqnarray*}
&&+\frac{n^{1/2}}{2}\left\{\mathbb{U}_n^{-1}\left(\frac{u k_n}{n}\right)-\frac{u k_n}{n}\right\}^{2}\frac{\partial^{2} \mathbb{C}}{\partial u^{2}}\left(u^{\prime},v^{\prime}\right)+\frac{n^{1/2}}{2}\left\{\mathbb{V}_n^{-1}\left(\frac{v k_n}{n}\right)-\frac{v k_n}{n}\right\}^{2}\frac{\partial^{2} \mathbb{C}}{\partial v^{2}}\left(u^{\prime},v^{\prime}\right)\\
&&+n^{1/2}\left\{\mathbb{U}_n^{-1}\left(\frac{u k_n}{n}\right)-\frac{u k_n}{n}\right\}\left\{\mathbb{V}_n^{-1}\left(\frac{v k_n}{n}\right)-\frac{v k_n}{n}\right\}\frac{\partial^2{} \mathbb{C}}{\partial u\partial v}\left(u^{\prime},v^{\prime}\right),
\end{eqnarray*}
where $(u^{\prime},v^{\prime})$ is a point between $\left(\frac{u k_n}{n},\frac{v k_n}{n}\right)$ and $\left(\mathbb{U}_n^{-1}\left(\frac{u k_n}{n}\right),\mathbb{V}_n^{-1}\left(\frac{v k_n}{n}\right)\right)$, and
$$
 \widetilde{\alpha}_n(u,v)= n^{1/2}(\mathbb{T}_n(u,v)-\mathbb{C}(u,v)).
$$
By using the preceding steps and facts, one finds
\begin{eqnarray*}
\mathbb{G}_{n}^{**}(u,v)&=& \widetilde{\alpha}_n\left(\mathbb{U}_n^{-1}\left(\frac{u
k_n}{n}\right),\mathbb{V}_n^{-1}\left(\frac{vk_n}{n}\right)\right)\\&& + n^{1/2}\left\{\mathbb{U}_n^{-1}\left(\frac{u k_n}{n}\right)-\frac{u k_n}{n}\right\}\frac{\partial \mathbb{C}}{\partial u}\left(\frac{u k_n}{n},\frac{vk_n}{n}\right)\\
&&+ n^{1/2}\left\{\mathbb{V}_n^{-1}\left(\frac{v k_n}{n}\right)-\frac{v k_n}{n}\right\}\frac{\partial \mathbb{C}}{\partial v}\left(\frac{u k_n}{n},\frac{vk_n}{n}\right)\\&&+\mathcal{O}(n^{-3/2}k_{n}(\log_{2}n)).\\
\end{eqnarray*}
For  each $n\geq 1$ and
$0\leq u,v \leq 1,$ set
\begin{eqnarray*}
\nonumber
\alpha_{n;0}^{**}(u,v)&:=& \widetilde{\alpha}_{n;0}\left(\frac{uk_n}{n},\frac{vk_n}{n} \right)\\
&:=& \widetilde{\alpha}_n\left(\frac{uk_n}{n},\frac{vk_n}{n}\right)-
\frac{\partial \mathbb{C}}{\partial v}\left(\frac{u k_n}{n},\frac{vk_n}{n}\right)
\alpha_{n,\mathbb{V}}\left(v\frac{k_n}{n}\right)\\&&-\frac{\partial \mathbb{C}}{\partial u}\left(\frac{u k_n}{n},\frac{vk_n}{n}\right)\alpha_{n,\mathbb{U}}\left(u\frac{k_n}{n}\right).
\end{eqnarray*}
For each $n\geq 1$ and $0 \leq u,v \leq 1$, we have
\begin{eqnarray*}
 \nonumber\lefteqn{\mathbb{G}_{n}^{**}(u,v)-\alpha_{n;0}^{**}(u,v)}\\
\nonumber & =&\left[\widetilde{\alpha}_n\left(\mathbb{U}_n^{-1}\left(\frac{u
k_n}{n}\right),\mathbb{V}_n^{-1}\left(\frac{vk_n}{n}\right)\right)-\widetilde{\alpha}_n \left(u\frac{k_n}{n},v\frac{k_n}{n} \right) \right] \\
\nonumber &  &              +
\frac{\partial \mathbb{C}}{\partial u}\left(\frac{u k_n}{n},\frac{vk_n}{n}\right)\left[\beta_{n;\mathbb{U}}\left(u\frac{k_n}{n}\right)+\alpha_{n;\mathbb{U}}\left(u\frac{k_n}{n}\right)\right]\\&&\nonumber
                 +\frac{\partial \mathbb{C}}{\partial v}\left(\frac{u k_n}{n},\frac{vk_n}{n}\right) \left[\beta_{n;\mathbb{V}}\left(v\frac{k_n}{n}\right)
                 +\alpha_{n;\mathbb{V}}\left(v\frac{k_n}{n}\right)\right] \\
\nonumber && + \mathcal{O}(n^{-3/2}k_{n}(\log_{2}n))\\
 &:= &R_{n;0}^{*}(u,v) + \frac{\partial \mathbb{C}}{\partial u}\left(\frac{u k_n}{n},\frac{vk_n}{n}\right) R_{n;\mathbb{U}}(u)
    \\&&+ \frac{\partial \mathbb{C}}{\partial v}\left(\frac{u k_n}{n},\frac{vk_n}{n}\right)  R_{n;\mathbb{V}}(v)+\mathcal{O}(n^{-3/2}k_{n}(\log_{2}n)).
\end{eqnarray*}
Using the fact that the first-order partial derivatives
of a copula are bounded (see for instance \cite{nelsen2006}), i.e.,
$$
0\leq  \frac{\partial \mathbb{C}(\cdot,\cdot)}{\partial u}\leq 1 ~~\mbox{and}~~0\leq  \frac{\partial \mathbb{C}(\cdot,\cdot)}{\partial v}\leq 1,
$$
 one can find the following
 \begin{eqnarray*}
 |\mathbb{G}_{n}^{**}(u,v)-\alpha_{n;0}^{**}(u,v)|
&\leq &|R_{n;0}^{*}(u,v)| + | R_{n;\mathbb{U}}(u)|+
|R_{n;\mathbb{V}}(v)|\\&&+\mathcal{O}(n^{-3/2}k_{n}(\log_{2}n)).
\end{eqnarray*}
Note that $|R_{n;0}^{*}(u,v)|$ may be treated making use the properties of the oscillation of the multivariate empirical process \cite[Theorem 1.7]{stute1984}.
By combining all this with the preceding proof, we obtain the following result for the general case of copulas.
\begin{theorem}\label{theoreme pavegene}
Let $\left\{k_n\right\}_{n\geq1}$ be a sequence of positive numbers fulfilling the assumptions
{\rm(H.1)-(H.4)}. Assume that $\mathbb{C}(\cdot,\cdot)$ is twice continuously differentiable on $(0,1)^2$ and all the
partial derivatives of second order are continuous on $[0,1]^2$. We have almost surely,
\begin{description}
    \item(i) when $ 0 < \gamma \leq 1/2$,
\begin{eqnarray} \label{resultat(i)G}
\nonumber&&\limsup_{n\rightarrow\infty} n^{1/2} k_n^{-1/4}(\log_2 n)^{-1/4}
(\log_{1} n)^{-1/2} \|\mathbb{G}_{n}^{**}-\alpha_{n;0}^{**} \|\\&&\qquad \qquad \leq
[3\times 2^{-1/4}+ \gamma 2^{5/4}](1-\gamma)^{1/4},
\end{eqnarray}
   \item(ii) when $ 1/2 < \gamma \leq 1$,
  \begin{eqnarray}\label{resultat(ii)G}
\nonumber&&\limsup_{n\rightarrow\infty} n^{1/2} k_n^{-1/4}(\log_2 n)^{-1/4}
(\log_{1} n)^{-1/2} \|\mathbb{G}_{n}^{**}-\alpha_{n;0}^{**} \|\\&&\qquad\qquad \leq
[3\times 2^{-3/4}+ \gamma 2^{3/4}] \gamma^{-1/4}.
\end{eqnarray}
\end{description}
\end{theorem}

 \begin{remark}
 In the general case, i.e., $\mathbb{C}(u,v)\neq uv$, the almost sure approximation rate is given, in \cite{Borisov1982}, by
 \begin{equation}\label{generalcase}
 \|\widetilde{\alpha}_{n}-\mathbf{B}_{n}^{*}\|=\mathcal{O}(n^{-1/6}\log n),
 \end{equation}
 where $\left\{\mathbf{B}_{n}^{*}(u,v): (u,v) \in[0,1]^{2}, n\geq 1\right\},$ is a  sequence of Brownian bridge, with covariance function
 $$ \mathbb{E}(\mathbf{B}_{n}^{*}(s_1,t_1)\mathbf{B}_{n}^{*}(s_2,t_2))= \mathbb{C}(s_1\wedge s_2,t_1\wedge t_2)-\mathbb{C}(s_1,t_1)\mathbb{C}(s_2,t_2) ,~~~~\hbox{for}~~(s_1,s_2,t_1,t_2)\in [0,1]^4.      $$
Note that the rate in (\ref{generalcase}) is not
$$
o\left(\frac{(\log n)^{1/2}(\log\log n)^{1/4}}{n^{1/4}}\right),
 $$
 then, for the empirical copula process, we have the following
 \begin{equation*}
 \|\mathbb{G}_{n}-\mathbb{B}_{n}^{*}\|=\mathcal{O}(n^{-1/6}\log n),
 \end{equation*}
 where
% \begin{equation}
%\widetilde{\mathbb{G}}_n(u,v):= n^{1/2}(\mathbb{C}_n(u,v)-
%\mathbb{C}(u,v)),~~~~\hbox{for}~~(u,v)\in [0,1]^2,
%\end{equation}
%and
 $$
 \mathbb{B}_{n}^{*}(s,t)=
\mathbf{B}_{n}^{*}(s,t)-\mathbf{B}_{n}^{*}(1,t)\frac{\partial \mathbb{C}}{\partial t}(s,t)-\mathbf{B}_{n}^{*}(s,1)\frac{\partial \mathbb{C}}{\partial s}(s ,t),~~~~\hbox{for}~~(s,t)\in [0,1]^2.
 $$
The results of Theorem \ref{theoreme pave} and Corollary \ref{COOOOO} are the best possible and are governed by the almost sure rate of Bahadur-Kiefer representation given in Fact 2.
\end{remark}

\subsection{ Application to tests of tail independence}\label{testing}
This section is largely inspired from \cite{deheuvels2006} and changes have been made in order to adopt it to our case.
Our main concern, is testing the null hypothesis $$\mathscr{H}_{0}~:~\mathbb{C}(u,v)=uv,~~\mbox{ for }~~	(u,v)\in[0, k_{n}/n]^{2}.$$
First, we introduce weighted bivariate tests of tail independence. Namely, for selected constants $\nu_{1}$ and $\nu_{2} \in \mathbb{R}$, we set
\begin{eqnarray}
\Omega_{n,k_{n},\nu_{1},\nu_{2}}=n\int_{0}^{k_{n}/n}\int_{0}^{k_{n}/n}u^{2\nu_{1}}v^{2\nu_{2}}\left\{\mathbb{C}_{n}(u,v)-uv\right\}^{2}dudv.
\end{eqnarray}
Recall the definition (\ref{GaussIndependen}) of $\mathbb{B}_{n}(\cdot,\cdot)$, for $n\geq 1$. Therefore, by the triangle inequality,
\begin{eqnarray}\label{triangl1}
\nonumber\|\mathbb{G}_n^{*2}-\mathbb{B}_{n}^{2}\| \leq \|\mathbb{G}_n^{*}-\mathbb{B}_{n}\|\times \{2\|\mathbb{G}_n^{*}\|+ \|\mathbb{G}_n^{*}-\mathbb{B}_{n}\|\}.
\end{eqnarray}
Note that the following elementary observation holds
\begin{equation}\label{weighted}
\int_{0}^{1} \int_{0}^{1} u^{2\nu_{1}}v^{2\nu_{2}}dudv<\infty,
\end{equation}
when $\nu_{1}>-1/2$ and $\nu_{2}>-1/2$. \\
We recall the following result given in \cite[Corollary 2.1]{deheuvels2006}. We have, with probability one,
\begin{equation}\label{lilG}
\limsup_{n\rightarrow \infty}(2\log\log n)^{-1/2}\|\mathbb{G}_{n}^{*}\|=\frac{1}{4}.
\end{equation}
By combining (\ref{triangl1}), (\ref{weighted}) and (\ref{lilG}), one finds the following.
\begin{corollary}
We have, almost surely, for $\nu_{1}>-1/2$ and $\nu_{2}>-1/2$,
\begin{description}
    \item(i)  when $ 0 < \gamma \leq 1/2$,
\begin{eqnarray} \label{resultat(i)}
\nonumber&&\limsup_{n\rightarrow\infty} n^{1/2} k_n^{-1/4}(\log_2 n)^{-3/4}
(\log_{1} n)^{-1/2} \left|\Omega_{n,k_{n},\nu_{1},\nu_{2}}-\int_{0}^{k_{n}/n}\int_{0}^{k_{n}/n}\mathbb{B}_{n}^{2}(u,v)dudv \right| \\&&\qquad\leq
[3\times 2^{-9/4}+ \gamma 2^{-3/4}](1-\gamma)^{1/4},
\end{eqnarray}
   \item(ii) when $ 1/2 < \gamma \leq 1$,
\begin{eqnarray} \label{resultat(ii)}
\nonumber&&\limsup_{n\rightarrow\infty} n^{1/2} k_n^{-1/4}(\log_2 n)^{-3/4}
(\log_{1} n)^{-1/2}  \left|\Omega_{n,k_{n},\nu_{1},\nu_{2}}-\int_{0}^{k_{n}/n}\int_{0}^{k_{n}/n}\mathbb{B}_{n}^{2}(u,v)dudv \right| \\&&\qquad\leq
[3\times 2^{-11/4}+ \gamma 2^{-5/4}] \gamma^{-1/4}.
\end{eqnarray}
\end{description}
\end{corollary}
\section{Strong approximation of smoothed local empirical copula process}\label{section3}
The smoothed empirical distribution function $\widehat{\mathbb{F}}_{n}(\cdot,\cdot)$ is defined,  for each $n\geq 1$ and $x,y \in \mathbb{R},$ by
$$
\widehat{\mathbb{F}}_{n}(x,y):=\frac{1}{n}\sum_{i=1}^{d}K_{n}(x-X_{i},y-Y_{i}).
$$
Here $K_{n}(x,y)=K(a_{n}^{-1/2}x,a_{n}^{-1/2}y)$, and
$$
K(x,y)=\int_{-\infty}^{x}\int_{-\infty}^{y}k(u,v)dudv
$$
for some bivariate kernel function $k:\mathbb{R}^{2}\mapsto\mathbb{R}  $, with $
\int_{-\infty}^{\infty}\int_{-\infty}^{\infty} k(x, y) dx dy= 1$, and a sequence of
bandwidths $a_{n}\downarrow 0$ as $n\rightarrow \infty$. For notational convenience, we have chosen the same
bandwidth sequence for each margin. This assumption can easily be dropped. For small
enough bandwidths $a_{n}$, the empirical distribution function $\mathbb{F}_{n}(\cdot,\cdot)$ and the smoothed empirical distribution function $\widetilde{\mathbb{F}}_{n}(\cdot,\cdot)$  are almost
indistinguishable, for more details refer to \cite[Lemma 7]{fermanianradulovicdragan2004} and also to \cite{vanderVaart1994}.
The continuity of the marginals $F(\cdot)$ and $G(\cdot)$ entails that we can replace them by uniform distributions. We then have, for each $n\geq 1$ and $0\leq u,v \leq 1$,
$$
\widehat{\mathbb{T}}_{n}(u,v):=\frac{1}{n}\sum_{i=1}^{d}K_{n}(u-U_{i},v-V_{i}),
$$
and the smoothed empirical copula function
$$
\widehat{\mathbb{C}}_n(u,v)=\widehat{\mathbb{T}}_{n}(G_{n}^{-1}(u),H_{n}^{-1}(v)).
$$
Consider the empirical process defined, for each $n\geq 1$ and $0\leq u,v \leq 1,$ by \begin{eqnarray}
\widehat{\alpha}_n(u,v)&=& n^{1/2}(\widehat{\mathbb{T}}_n(u,v)-uv),\\
\widehat{\alpha}_n^{*}(u,v)&:=& \widehat{\alpha}_n\left(u\frac{k_n}{n},v\frac{k_n}{n}\right),
\end{eqnarray}
and  define the smoothed empirical copula process
 $\widehat{\mathbb{G}}_n(\cdot,\cdot)$ by setting
\begin{equation}\label{processus de copule}
\widehat{\mathbb{G}}_n(u,v):= n^{1/2}(\widehat{\mathbb{C}}_n(u,v)-
uv),~~~~\hbox{for}~~(u,v)\in [0,1]^2.
\end{equation}
The main aim of this section is to investigate the smoothed empirical copula process defined, in terms of a sequence $\{k_n\}_{n=1}^{\infty}$, for each $n\geq 1$,  by
\begin{eqnarray} \label{processus G_n lisse sur le pave}
\widehat{\mathbb{G}}_n^{*}(u,v):=
\widehat{\mathbb{G}}_n\left(u\frac{k_n}{n},v\frac{k_n}{n}\right)~~~~\hbox{for}~~0\leq
u,v \leq 1.
\end{eqnarray}
Observe that
\begin{eqnarray}
\nonumber
\widehat{\mathbb{G}}_{n}^{*}(u,v)&=& n^{1/2}\left\{\widehat{\mathbb{T}}_n\left(\mathbb{U}_n^{-1}\left(\frac{u k_n}{n}\right),\mathbb{V}_n^{-1}\left(\frac{vk_n}{n}\right)\right)-\frac{u k_n}{n}\frac{vk_n}{n} \right\}\\
\nonumber &=& \widehat{\alpha}_n\left(\mathbb{U}_n^{-1}\left(\frac{u
k_n}{n}\right),\mathbb{V}_n^{-1}\left(\frac{vk_n}{n}\right)\right)\\&& \nonumber+ n^{1/2}\left[\mathbb{U}_n^{-1}\left(\frac{u k_n}{n}\right)\mathbb{V}_n^{-1}\left(\frac{vk_n}{n}\right)-\frac{u k_n}{n}\frac{vk_n}{n}\right]\\
\nonumber &=& \widehat{\alpha}_n\left(\mathbb{U}_n^{-1}\left(\frac{u
k_n}{n}\right),\mathbb{V}_n^{-1}\left(\frac{vk_n}{n}\right)\right)+u
\frac{k_n}{n} \beta_{n;\mathbb{V}}\left(v\frac{k_n}{n}\right)\\&&+v\frac{k_n}{n}\beta_{n;\mathbb{U}}\left(u\frac{k_n}{n}\right)+n^{-1/2}\beta_{n;\mathbb{U}}\left(u\frac{k_n}{n}\right)\beta_{n;\mathbb{V}}\left(v\frac{k_n}{n}\right).
\end{eqnarray}
We shall assume that
\begin{equation}\label{conditionK}
\int_{-\infty}^{\infty}\int_{-\infty}^{\infty}(x^{2}+y^{2})^{1/2}dK(x,y)<\infty.
\end{equation}
Since $\mathbb{T}(\cdot,\cdot)$ satisfies condition (2.1) in \cite{Mason2010}, one can apply \cite[Corollary 2, eq (2.7)]{Mason2010} to obtain
\begin{equation}
\|\widehat{\mathbb{T}}_{n} -\mathbb{T}_{n}\|=O\left(\frac{\sqrt{a_{n}\log_{2} n}}{\sqrt{n}}\right),~~~a.s.,
\end{equation}
which gives
\begin{eqnarray}
\nonumber
\widehat{\mathbb{G}}_{n}^{*}(u,v) &=& \alpha_n\left(\mathbb{U}_n^{-1}\left(\frac{u
k_n}{n}\right),\mathbb{V}_n^{-1}\left(\frac{vk_n}{n}\right)\right)\\&&+u
\frac{k_n}{n} \beta_{n;\mathbb{V}}\left(v\frac{k_n}{n}\right)+v\frac{k_n}{n}\beta_{n;\mathbb{U}}\left(u\frac{k_n}{n}\right)\nonumber\\
&&\nonumber+n^{-1/2}\beta_{n;\mathbb{U}}\left(u\frac{k_n}{n}\right)\beta_{n;\mathbb{V}}\left(v\frac{k_n}{n}\right)\\&&+O\left(\sqrt{a_{n}\log_{2} n}\right),~~~a.s.
\end{eqnarray}
Here is our main result concerning  $\widehat{\mathbb{G}}_{n}^{*}(\cdot,\cdot)$.
\begin{theorem}\label{theoreme pavelisse}
Let $\left\{k_n\right\}_{n\geq1}$ be a sequence of positive numbers fulfilling the assumptions
{\rm(H.1)-(H.4)} and $K(\cdot,\cdot)$ satisfies (\ref{conditionK}). For $a_{n}=O(n^{-(1/4+\delta)}),$ $\delta>0$, we have almost surely,
\begin{description}
    \item(i)  when $ 0 < \gamma \leq 1/2$,
\begin{eqnarray}
\nonumber&&\limsup_{n\rightarrow\infty} n^{1/2} k_n^{-1/4}(\log_2 n)^{-1/4}
(\log_{1}n)^{-1/2} \|\widehat{\mathbb{G}}_{n}^{*}-\alpha_{n;0}^{*} \|\\&&\qquad\qquad \leq
[3\times 2^{-1/4}+ \gamma 2^{5/4}](1-\gamma)^{1/4},
\end{eqnarray}
   \item(ii) when $ 1/2 < \gamma \leq 1$,
  \begin{eqnarray}
\nonumber&&\limsup_{n\rightarrow\infty} n^{1/2} k_n^{-1/4}(\log_2 n)^{-1/4}
(\log_{1}n)^{-1/2} \|\widehat{\mathbb{G}}_{n}^{*}-\alpha_{n;0}^{*} \|\\ &&\qquad\qquad\leq
[3\times 2^{-3/4}+ \gamma 2^{3/4}] \gamma^{-1/4}.
\end{eqnarray}
\end{description}
\end{theorem}
We have immediately.
\begin{corollary} Under the conditions of preceding theorem and Theorem \ref{bestapproximation}, we have almost surely,
\begin{description}
    \item(i)  when $ 0 < \gamma \leq 1/2$,
\begin{eqnarray}
\nonumber&&\limsup_{n\rightarrow\infty} n^{1/2} k_n^{-1/4}(\log_2 n)^{-1/4}
(\log_{1} n)^{-1/2} \|\widehat{\mathbb{G}}_{n}^{*}-\mathbb{B}_{n} \| \\&&\qquad\qquad\leq
[3\times 2^{-1/4}+ \gamma 2^{5/4}](1-\gamma)^{1/4},
\end{eqnarray}
   \item(ii) when $ 1/2 < \gamma \leq 1$,
  \begin{eqnarray}
\nonumber&&\limsup_{n\rightarrow\infty} n^{1/2} k_n^{-1/4}(\log_2 n)^{-1/4}
(\log_{1} n)^{-1/2} \|\widehat{\mathbb{G}}_{n}^{*}-\mathbb{B}_{n} \| \\&&\qquad\qquad\leq
[3\times 2^{-3/4}+ \gamma 2^{3/4}] \gamma^{-1/4}.
\end{eqnarray}
\end{description}
\end{corollary}

%\section*{Aknowledgment}
%We are grateful for the comments and suggestions made by a referee that helped to improve the presentation of our results.


\begin{thebibliography}{}
\bibitem[Borisov(1982)]{Borisov1982}
Borisov, I.~S. (1982).
\newblock Approximation of empirical fields that are constructed from vector
  observations with dependent coordinates.
\newblock {\em Sibirsk. Mat. Zh.}, {\bf 23}(5), 31--41, 222.

\bibitem[Castelle and Laurent-Bonvalot(1998)]{castellebonvalot1998}
Castelle, N. and Laurent-Bonvalot, F. (1998).
\newblock Strong approximations of bivariate uniform empirical processes.
\newblock {\em Ann. Inst. H. Poincar\'e Probab. Statist.}, {\bf 34}(4),
  425--480.

\bibitem[Cs{\"o}rg{\H{o}} and Horv{\'a}th(1988)]{Csorgo1988}
Cs{\"o}rg{\H{o}}, M. and Horv{\'a}th, L. (1988).
\newblock Asymptotic tail behavior of uniform multivariate empirical processes.
\newblock {\em Ann. Prob.}, {\bf 18}(4), 1723--1738.

\bibitem[Deheuvels(1979)]{deheuvels1979a}
Deheuvels, P. (1979).
\newblock La fonction de d\'ependance empirique et ses propri\'et\'es. {U}n
  test non param\'etrique d'ind\'ependance.
\newblock {\em Acad. Roy. Belg. Bull. Cl. Sci. (5)}, {\bf 65}(6), 274--292.

\bibitem[Deheuvels(1980)]{Deheuvels1981b}
Deheuvels, P. (1980).
\newblock Nonparametric test of independence.
\newblock In {\em Nonparametric asymptotic statistics (Proc. Conf., Rouen,
  1979) (French)}, volume 821 of {\em Lecture Notes in Math.}, pages 95--107.
  Springer, Berlin.

\bibitem[Deheuvels(1981)]{deheuvels1981}
Deheuvels, P. (1981).
\newblock Multivariate tests of independence.
\newblock In {\em Analytical methods in probability theory (Oberwolfach,
  1980)}, volume 861 of {\em Lecture Notes in Math.}, pages 42--50. Springer,
  Berlin.

\bibitem[Deheuvels(1997)]{sDeheuvels1997}
Deheuvels, P. (1997).
\newblock Strong laws for local quantile processes.
\newblock {\em Ann. Probab.}, {\bf 25}(4), 2007--2054.

\bibitem[Deheuvels {\em et~al.}(2006)]{deheuvels2006}
Deheuvels, P., Peccati, G., and Yor, M. (2006).
\newblock On quadratic functionals of the {B}rownian sheet and related
  processes.
\newblock {\em Stochastic Process. Appl.}, {\bf 116}(3), 493--538.

\bibitem[Deheuvels(2007)]{Deheuvels2007}
Deheuvels, P. (2007).
\newblock Weighted multivariate tests of independence.
\newblock {\em Comm. Statist. Theory Methods}, {\bf 36}(13-16), 2477--2491.


\bibitem[Deheuvels(2009)]{Deheuvelsref2009}
Deheuvels, P. (2009).
\newblock A multivariate Bahadur-Kiefer representation for the empirical copula
  process.
\newblock {\em Zap. Nauchn. Sem. S.-Peterburg. Otdel. Mat. Inst. Steklov.
  (POMI)}, {\bf 364}, 120--147.

\bibitem[Einmahl and Mason(1988a)]{masoneinmahl1988}
Einmahl, J. H.~J. and Mason, D.~M. (1988a).
\newblock Laws of the iterated logarithm in the tails for weighted uniform
  empirical processes.
\newblock {\em Ann. Probab.}, {\bf 16}(1), 126--141.

\bibitem[Einmahl and Mason(1988b)]{masonein1988}
Einmahl, J. H.~J. and Mason, D.~M. (1988b).
\newblock Strong limit theorems for weighted quantile processes.
\newblock {\em Ann. Probab.}, {\bf 16}(4), 1623--1643.


\bibitem[Einmahl and Ruymgaart(1987)]{sEinmahl1987}
Einmahl, J. H.~J. and Ruymgaart, F.~H. (1987).
\newblock The almost sure behavior of the oscillation modulus of the
  multivariate empirical process.
\newblock {\em Statist. Probab. Lett.}, {\bf 6}(2), 87--96.

\bibitem[Fermanian {\em et~al.}(2004)]{fermanianradulovicdragan2004}
Fermanian, J.-D., Radulovi{\'c}, D., and Wegkamp, M. (2004).
\newblock Weak convergence of empirical copula processes.
\newblock {\em Bernoulli}, {\bf 10}(5), 847--860.

\bibitem[Joe(1997)]{Joe1997}
Joe, H. (1997).
\newblock {\em Multivariate models and dependence concepts}, volume~73 of {\em
  Monographs on Statistics and Applied Probability}.
\newblock Chapman \& Hall, London.



\bibitem[Gaenssler and Stute(1987)]{stute1987}
Gaenssler, P. and Stute, W. (1987).
\newblock {\em Seminar on empirical processes}, volume~9 of {\em DMV Seminar}.
\newblock Birkh\"auser Verlag, Basel.

\bibitem[Mason and Swanepoel(2010)]{Mason2010}
Mason, D.~M. and Swanepoel, J. W.~H. (2010).
\newblock A general result on the uniform in bandwidth consistency of
  kernel-type function estimators.
\newblock {\em TEST}, {\bf 19}(2), DOI: 10.1007/s11749-010-0188-0.

\bibitem[Mason and van Zwet(1987)]{Mason1987}
Mason, D.~M. and van Zwet, W.~R. (1987).
\newblock A refinement of the {KMT} inequality for the uniform empirical
  process.
\newblock {\em Ann. Probab.}, {\bf 15}(3), 871--884.

\bibitem[Mason (1988)]{Mason1988}
Mason, D.~M.(1988).
\newblock A strong invariance theorem for the tail empirical process.
\newblock {\em Ann. Inst. Henri Poincar\'e}, {\bf 24}(4), 491--506.

\bibitem[Moore and Spruill(1975)]{Moore1975}
Moore, D.~S. and Spruill, M.~C. (1975).
\newblock Unified large-sample theory of general chi-squared statistics for
  tests of fit.
\newblock {\em Ann. Statist.}, {\bf 3}, 599--616.


\bibitem[Nelsen(2006)]{nelsen2006}
Nelsen, R.~B. (2006).
\newblock {\em An introduction to copulas}.
\newblock Springer Series in Statistics. Springer, New York, second edition.

\bibitem[R{\"u}schendorf(1974)]{Ruch1}
R{\"u}schendorf, L. (1974).
\newblock On the empirical process of multivariate, dependent random variables.
\newblock {\em J. Multivariate Anal.}, {\bf 4}, 469--478.

\bibitem[R{\"u}schendorf(1976)]{Rush2}
R{\"u}schendorf, L. (1976).
\newblock Asymptotic distributions of multivariate rank order statistics.
\newblock {\em Ann. Statist.}, {\bf 4}(5), 912--923.

\bibitem[Sklar(1959)]{Sklar1959}
Sklar, M. (1959).
\newblock Fonctions de r\'epartition \`a {$n$} dimensions et leurs marges.
\newblock {\em Publ. Inst. Statist. Univ. Paris}, {\bf 8}, 229--231.

\bibitem[Sklar(1973)]{Sklar1973}
Sklar, A. (1973).
\newblock Random variables, joint distribution functions, and copulas.
\newblock {\em Kybernetika (Prague)}, {\bf 9}, 449--460.


\bibitem[Stute(1984)]{stute1984}
Stute, W. (1984).
\newblock {The oscillation behavior of empirical processes: The multivariate
  case.}
\newblock {\em Ann. Probab.}, {\bf 12}, 361--379.



\bibitem[van~der Vaart(1994)]{vanderVaart1994}
van~der Vaart, A. (1994).
\newblock Weak convergence of smoothed empirical processes.
\newblock {\em Scand. J. Statist.}, {\bf 21}(4), 501--504.

\bibitem[Van~der Vaart and Wellner(1996)]{Wellner1996}
van~der Vaart, A.~W. and Wellner, J.~A. (1996).
\newblock {\em Weak convergence and empirical processes}.
\newblock Springer Series in Statistics. Springer-Verlag, New York.
\newblock With applications to statistics.


\end{thebibliography}
\end{document}